\documentclass{svmult}
\usepackage{amssymb}
\usepackage{amsmath}
\usepackage{graphicx}

\DeclareMathOperator{\zr}{\mathbb{Z}}
\DeclareMathOperator{\re}{\mathbb{R}}
\DeclareMathOperator{\p}{\mathbb{P}}
\DeclareMathOperator{\e}{\mathbb{E}}

\newtheorem{prop}[theorem]{Proposition}

\date{January, 2005}
\begin{document}
\title*{The maximality principle revisited: on certain optimal stopping problems}
\author{Jan Ob\l \'oj\thanks{ Work partially supported by Polish Academy of Science scholarship.}}
\titlerunning{On certain optimal stopping problems}
\institute{Laboratoire de Probabilit\'es et Mod\`eles Al\'eatoires, Universit\'e Paris 6,\newline\hspace*{0.7cm}4 pl. Jussieu, Bo\^{i}te 188, 75252 Paris Cedex 05, France\newline Department of Mathematics, Warsaw University\newline\hspace*{0.7cm}Banacha 2, 02-097 Warszawa, Poland; \texttt{obloj@mimuw.edu.pl}}
\maketitle
\abstract{We investigate in detail works of Peskir  \cite{MR2000c:60052} and Meilijson  \cite{optimal_meilijson} and develop a link between them.
We consider the following optimal stopping problem: maximize $V_\tau=\e\Big[\phi(S_\tau)-\int_0^\tau c(B_s)ds\Big]$ over all stopping times with $\e\int_0^\tau c(B_s)ds<\infty$, where $S=(S_t)_{t\geq 0}$ is the maximum process associated with real valued Brownian motion $B$, $\phi\in C^1$ is non-decreasing and $c\ge 0$ is continuous. From work of Peskir \cite{MR2000c:60052} we deduce that this problem has a unique solution if and only if the differential equation
$$g'(s)=\frac{\phi'(s)}{2c(g(s))(s-g(s))}$$
admits a maximal solution $g_*(s)$ such that $g_*(s)< s$ for all $s\geq 0$. The stopping time which yields the highest payoff
 can be written as $\tau_*=\inf\{t\geq 0: B_t\leq g_*(S_t)\}$. The problem is actually solved in a general case of a real-valued, time homogeneous diffusion $X=(X_t:t\geq 0)$ instead of $B$.
We then proceed to solve the problem for more general functions $\phi$ and $c$.
Explicit formulae for payoff are given.

We apply the results to solve the so-called optimal Skorokhod embedding problem. We give also a sample of applications to various inequalities dealing with terminal value and maximum of a process.
}
\newpage

\section{Introduction}
In this paper we study certain optimal stopping problems. Our interest in these problems originated from the Skorokhod embedding field and a certain knowledge of the latter will be useful for reading (see our survey paper \cite{genealogia}).

We are directly stimulated by an article by Meilijson \cite{optimal_meilijson} and our observation that his results can be seen via the tools developed by Dubins, Shepp and Shiryaev \cite{MR96j:60077} and Peskir \cite{MR2000c:60052}. In this paper, we manage to establish a link and generalize both the works of Meilijson \cite{optimal_meilijson} and Peskir \cite{MR2000c:60052}. We try to investigate the nature of the optimal stopping problem (\ref{eq:optimal_stopping}) through a series of remarks and rely on reasonings found in the above-cited articles together with some approximation techniques and limit passages to develop certain generalizations. Our work has therefore characteristics of an exposition of the subject with aim at unifying and extending known results.

Let $\phi$ be a non-negative, increasing, continuous function and $c$ a continuous, positive function, and consider the following optimal stopping problem of maximizing
\begin{equation}
  \label{eq:optimal_stopping1}
  V_\tau=\e\Big[\phi(S_\tau)-\int_0^\tau c(B_s)ds\Big],
\end{equation}
over all stopping times $\tau$ such that
\begin{equation}
  \label{eq:opt_stop_cond}
  \e\Big[\int_0^\tau c(B_s)ds\Big]<\infty,
\end{equation}
where $(B_t:t\ge 0)$ is a real-valued Brownian motion and $S_t$ is its unilateral maximum, $S_t=\sup_{u\le t}B_y$.

Suppose, in the first moment, that $\phi(x)=x$ and $c(x)=c>0$ is a constant. In this formulation the problem was solved by Dubins and Schwarz \cite{MR89m:60101} in an article on Doob-like inequalities. The optimal stopping time is just the Az\'ema-Yor embedding (see Az\'ema and Yor \cite{MR82c:60073a} or Section 5 in Ob\l \'oj \cite{genealogia}) for a shifted (to be centered) exponential distribution with parameter $2c$. This leads in particular to an optimal inequality (\ref{eq:ineq1}).

Keeping $\phi(x)=x$, let $c(x)$ be a non-negative, continuous function.
The setup was treated by Dubins, Shepp and Shiryaev \cite{MR96j:60077}, and by Peskir in a series of articles (\cite{MR2000c:60052}, \cite{MR2001e:60083}, \cite{peskir_maxprocinopt}). Peskir treated the case of real diffusions, which allows to recover the solution for general $\phi$ as a corollary.
\begin{theorem}[Peskir \cite{MR2000c:60052}]
\label{thm:maximality_principle}
The problem of maximizing (\ref{eq:optimal_stopping1}) over all stopping times $\tau$ satisfying (\ref{eq:opt_stop_cond}), for $\phi(x)=x$ and $c(x)$ a non-negative, continuous function, has an optimal solution with finite payoff if and only if there exists a maximal solution $g_*$ of
\begin{equation}
  \label{eq:diff_g}
  g'(s)=\frac{1}{2c(g(s))(s-g(s))}
\end{equation}
  which stays strictly below the diagonal in $\re^2$. The Az\'ema-Yor stopping time $\tau_*=\inf\{t\geq 0: B_t\leq g_*(S_t)\}$ is then optimal and satisfies (\ref{eq:opt_stop_cond}) whenever there exists a stopping time which satisfies (\ref{eq:opt_stop_cond}).
\end{theorem}
As pointed out above, this theorem was proved in a setup of any real, regular, time-homogeneous diffusion to which we will come back later. The characterization of existence of a solution to (\ref{eq:optimal_stopping1}) through existence of a solution to the differential equation (\ref{eq:diff_g}) is called the \textit{maximality principle}. We point out that Dubins, Shepp and Shiryaev, who worked with Bessel processes, had a different way of characterizing the optimal solution to (\ref{eq:diff_g}), namely they required that $\frac{g_*(s)}{s}\xrightarrow[s\to\infty]{} 1$.

Let now $\phi$ be any non-negative, non-decreasing, right-continuous function such that $\phi(B_t)-ct$ is a.s. negative on some $(t_0,\infty)$ (with $t_0$ random) and keep $c$ constant. This optimal stopping problem was solved by Meilijson \cite{optimal_meilijson}. Define
$$H(x)=\sup_\tau\e\Big[\phi(x+S_\tau)-c\tau\Big].$$
\begin{theorem}[Meilijson \cite{optimal_meilijson}]
\label{thm:meilijson1}
  Suppose that $\e\sup_t\{\phi(B_t)-ct\}<\infty$. Then $H$ is absolutely continuous and is the minimal solution to the differential equation
  \begin{equation}
    \label{eq:diff_H}
    H(x)-\frac{1}{4c}(H'(x))^2=\phi(x)\;.
  \end{equation}
If $\phi$ is constant on $[x_0,\infty)$ then $H$ is the unique solution to (\ref{eq:diff_H}) that equals $\phi$ on $[x_0,\infty)$. The optimal stopping time $\tau_*$ which yields $H(0)$ is the Az\'ema-Yor stopping time given by $\tau_*=\inf\{t\geq 0: B_t\leq S_t-\frac{H'(S_t)}{2c}\}$.
\end{theorem}

Let us examine in more detail the result of Meilijson in order to compare it with the result of Peskir. The Az\'ema-Yor stopping time is defined as $\tau_*=\inf\{t\geq 0: B_t\leq g(S_t)\}$ with $g(x)=x-\frac{H'(x)}{2c}$. Let us determine the differential equation satisfied by $g$. Note that $H$ is by definition non-decreasing, so we have $H'(x)=\sqrt{4c}\sqrt{(H(x)-\phi(x))}$. For suitable $\phi$, this is a differentiable function and differentiating it we obtain
\begin{eqnarray}
\label{eq:diff_H_bis}
  H''(x)&=&\frac{2c(H'(x)-\phi'(x))}{\sqrt{4c(H(x)-\phi(x))}}=\frac{2c(H'(x)-\phi'(x))}{H'(x)}\nonumber.
\end{eqnarray}
Therefore
\begin{equation}
  \label{eq:diff_g_meilij}
  g'(x)=1-\frac{H''(x)}{2c}=\frac{\phi'(x)}{H'(x)}=\frac{\phi'(x)}{2c(x-g(x))}.
\end{equation}
We recognize immediately the equation (\ref{eq:diff_g}) only there $\phi'(s)=1$ and $c$ was a function and not a constant.
This motivated our investigation of the generalization of problem (\ref{eq:optimal_stopping1}), given in (\ref{eq:optimal_stopping}), which is solved in Theorems \ref{thm:gen_max1} - \ref{thm:gen_max3} in Section \ref{sec:gen_max}.

Consider now the converse problem. That is, given a centered probability measure $\mu$ describe all pairs of functions $(\phi,c)$ such that the optimal stopping time $\tau_*$, which solves (\ref{eq:optimal_stopping1}) exists and embeds $\mu$, that is $B_{\tau_\mu}\sim\mu$.

This was called the ``optimal Skorokhod embedding problem'' (term introduced by Peskir in \cite{MR2001e:60083}), as we not only specify a method to obtain an embedding for $\mu$ but also construct an optimal stopping problem, of which this embedding is a solution.

Again, let us consider two special cases. First, let $\phi(x)=x$.
From the Theorem \ref{thm:maximality_principle} we know that $\tau_*$ is the Az\'ema-Yor stopping time and the function $g_*$ is just the inverse of barycentre function of some measure $\mu$. The problem therefore consists in identifying the dependence between $\mu$ and $c$.
Suppose for simplicity that $\mu$ has a strictly positive density $f$ and that it satisfies $L\log L$-integrability condition.
\begin{theorem}[Peskir \cite{MR2001e:60083}]
\label{thm:peskir2}
In the setup above there exists a unique function
\begin{equation}
  \label{eq:hazard_function}
  c(x)=\frac{1}{2}\frac{f(x)}{\overline{\mu}(x)}
\end{equation}
such that the optimal solution $\tau_*$ of (\ref{eq:optimal_stopping1}) embeds $\mu$. This optimal solution is then the Az\'ema-Yor stopping time given by
\begin{equation}
  \label{eq:ay_st_def}
\tau_*=\inf\{t\geq 0:S_t\geq \Psi_\mu(B_t)\}\textrm{, where }
\end{equation}
\begin{equation}
  \label{eq:ay_def}
  \Psi_\mu(x)=\frac{1}{\mu([x,\infty))}\int_{[x,\infty)}yd\mu(y).
\end{equation}
\end{theorem}

\noindent The function $c$ in (\ref{eq:hazard_function}) is recognized as a half of the hazard function which plays an important role in some studies. Now let $c(x)=c$ be constant. Then we have
\begin{theorem}[Meilijson \cite{optimal_meilijson}]
In the setup above, there exists a unique function $\phi$ defined through (\ref{eq:diff_H}) with $H'(x)=2c(x-\Psi_\mu^{-1}(x))$, where $\Psi_\mu$ is the barycentre function given in (\ref{eq:ay_def}), such that
the optimal solution $\tau_*$ of (\ref{eq:optimal_stopping}) embeds $\mu$. This optimal solution is then the Az\'ema-Yor stopping time given by
(\ref{eq:ay_st_def}).
\end{theorem}
We will provide a general solution that identifies all the pairs $(\phi,c)$ with basic regularity properties, in Proposition \ref{prop:opt_skoro} in Section \ref{sec:gen_skoro}.

\section{Notation}
\label{sec:notation}

In this section we fix the notation that, unless stated otherwise, is used in the rest of the paper.

$X=(X_t)_{t\geq 0}$ denotes a real-valued, time-homogeneous, regular diffusion associated with the infinitesimal generator
\begin{equation}
  \label{eq:infi_gen}
  \mathbb{L}_X=\delta(x)\frac{\partial}{\partial x}+\frac{\sigma^2(x)}{2}\frac{\partial^2}{\partial x^2},
\end{equation}
where the drift coefficient $x\to \delta(x)$ and the diffusion coefficient $x\to \sigma(x)>0$ are continuous. We assume moreover that there exists a real-valued, standard Brownian motion $B=(B_t)_{t\geq 0}$, such that $X$ solves the stochastic differential equation
\begin{equation}
  \label{eq:diff_eq_X}
  d X_t=\delta(X_t)dt+\sigma(X_t)dB_t,
\end{equation}
with $X_0=x$ under $\p_x:=\p$ for $x\in \re$. The natural filtration of $X$ is denoted $(\mathcal{F}_t: t\geq 0)$ and is taken right-continuous and completed. The scale function and the speed measure of $X$ are denoted respectively by $L$ and by $m$. The state space of $X$ is $(a_X,b_X)$, where $a_X<b_X$ can be finite or infinite. The interval $(a_X,b_X)$ might be also closed from one, or both, sides under some additional conditions, as will be discussed later.
We recall, that for Brownian motion $L(x)=x$ and $m(dx)=2dx$.

The one-sided maximum process associated with $X$ is
\begin{equation}
  \label{eq:deff_S}
  S_t=\Big(\sup_{0\leq r\leq t}X_r\Big)\lor s
\end{equation}
started at $s\geq x$ under $\p_{x,s}:=\p$. When several maximum processes will appear for different diffusions, we will add appropriate superscripts, thus writing $S^X_t$. The first hitting times for $S$ are noted $T_a=\inf\{t\geq 0:S_t\geq a\}$.

We make the following assumptions on the function $\phi$:
\begin{itemize}
\item $\phi:(a_X,b_X)\to \re$ is non-decreasing, right-continuous and its points of discontinuity are isolated;
\item $\phi'$, taken right-continuous, is well defined and its points of discontinuity are isolated;
\item there exists $r_\phi$, such that on $[r_\phi-1,+\infty)$, $\phi\in C^1$ and $\phi'>0$ ;
\end{itemize}
The last assumption has technical signification and is there to allow for a convenient identification of a particular solution to a differential equation. It will be clear form our proofs of main theorems, we hope, that one could impose some other condition, which would yield a different description.

We denote $D_\phi$, the set of points of discontinuity of $\phi$. It is bounded from above by $r_\phi$ and can write $D_\phi=\{d_0,d_1,\dots\}$ where $d_i\geq d_{i+1}$. Similarly, the set of points of discontinuity of $\phi'$ is also bounded from above by $r_\phi$ and we can write $D_{\phi'}=D_\phi\cup \{f_0,f_1,\dots\}$, where $f_i\geq f_{i+1}$.

It will be convenient to work sometimes with a different function $\tilde{\phi}$, which is of class $C^1$ with $\tilde{\phi}'>0$, and satisfies: $\tilde{\phi}\leq \phi$, $\tilde{\phi}_{|_{[r_\phi,\infty)}}=\phi_{|_{[r_\phi,\infty)}}$. No reasoning will depend on a particular choice of $\tilde{\phi}$, so we do not present a specific construction.

The function $c:[0,+\infty)\to [0,+\infty)$ has a countable number of discontinuities.
In his work Peskir \cite{MR2000c:60052} supposed $c>0$, but we will see (cf. Remark \ref{rq:c}) that allowing $c$ to be zero on some intervals can be of great use. In some parts of the work we will make additional assumptions on the regularity of $c$, yet in others we will allow it to take the value $+\infty$.

We consider the following optimal stopping problem:
\begin{equation}
  \label{eq:optimal_stopping}
  V_*(x,s)=\sup_{\tau}\e_{x,s}\Big(\phi(S_\tau)-\int_0^\tau c(X_t)dt\Big),
\end{equation}
where $\tau$ is a $\mathcal{F}_t$-stopping time, which satisfies
\begin{equation}
    \label{eq:optimal_stopping_cond}
\e_{x,s}\Big(\int_0^\tau c(X_t)dt\Big)<\infty.
\end{equation}
Through a solution to ``the optimal stopping problem (\ref{eq:optimal_stopping})'' we understand a stopping time which yields $V_*$ and satisfies the condition (\ref{eq:optimal_stopping_cond}).
Such solution is denoted $\tau_*$, i.e. $V_*(x,s)=\e_{x,s}\Big(\phi(S_{\tau_*})-\int_0^{\tau_*} c(X_t)dt\Big)$. Note that of course $V_*$ and $\tau_*$ depend on $\phi$ and $c$ but it will be obvious from the context if we discuss the general setup of arbitrary $(\phi,c)$ or some special cases.

\section{Some remarks on the problem}
\label{sec:remarks}

In the next section we will present a complete solution to the optimal stopping problem (\ref{eq:optimal_stopping}). However, before we do it, we want to give some relatively simple observations, which help to understand better the nature of the problem. Most of them become nearly evident once a person is acquainted with the problem, yet we think it is worthwhile to gather them here. As a matter of fact, they will help us a great deal, both in formulating and in proving, the main results of this paper.

The general construction of the optimal stopping problem under consideration is such, that we get rewarded according to the record-high value so far, and get punished all the time, proportionally to the time elapsed and depending on the path of our process. This means basically that the situation when the diffusion increases is favorable and what is potentially dangerous, as it might be ``too costly'', are excursions far away below the maximum, that is the negative excursions of the process $(X_t-S_t)$. The first observation is then that we should not stop in the support of $dS_t$.

\begin{remark}[Proposition 2.1 \cite{MR2000c:60052}]
  If $\phi'>0$ and $\tau$ yields the solution to the problem (\ref{eq:optimal_stopping}), then $X_\tau\neq S_\tau$ a.s., that is the process $(X_t,S_t)$ cannot be optimally stopped on the diagonal of $\re^2$.
\end{remark}

Let us investigate in more detail the nature of the optimal stopping time. If we follow the process $X_t$, as noted above, we gain when it increases. In contrast, when it decreases we only get punished all the time. This means that at some point it gets too costly to continue a given negative excursion of $(X_t-S_t)$. Exactly when it becomes too costly to continue the excursion depends only on three factors: the diffusion characteristics, the functions $\phi$, and $c$, and on the value $S_t$ of the maximum so far. Indeed, the choice of the point in an obvious manner depends on what we expect to earn, if the process climbs back to $S_t$, which in turn is a function of the diffusion characteristics and the functions $\phi$, and $c$. Note however, that through the Markov property, it doesn't depend on the way the process arrives at $S_t$. Similarly, thanks to the Markov property, the way the excursion of $(X-S)$ straddling time $t$ develops, doesn't affect the choice of the stopping value. Clearly also this value is an increasing function of the maximum.
\begin{remark}
\label{rq:formofthesolution}
  The solution to the optimal stopping problem (\ref{eq:optimal_stopping}), if it exists, is given by a stopping time of the following form $\tau=\inf\{t\geq 0: X_t\leq g_*(S_t)\}$, where the function $g_*$ is non-decreasing and depends on the characteristics of $X$ and the functions $\phi$, and $c$.
\end{remark}
Now that we have a feeling of how the stopping time looks like, we will try to learn more about the function $g_*$.
Consider the function $V(x,s)$ given by (\ref{eq:optimal_stopping}). Since $V(x,s)$ corresponds to the optimal stopping of the process $(X_t,S_t)$ started at $(x,s)$, the values of $\phi(u)$ for $u<s$ never intervene. Furthermore, if the process reaches some value $s_1>s$ then, thanks to the strong Markov property, we just face the problem $V(s_1,s_1)$. We could then replace the function $\phi$ with some other function $\tilde{\phi}$ which coincides with $\phi$ on $[s_1,\infty)$ and this would not affect the stopping rule, that is it would not affect values of the function $g_*$ on the interval $[s_1,\infty)$. Thus the finiteness of $V$ is not affected by any change of $\phi$ away from infinity. Likewise, if we add a constant $\delta$ to $\phi$ it will only change $V(\cdot,\cdot)$ by $\delta$ and leave $g_*$ unchanged. Combining these observations we arrive at the following remark.
\begin{remark}
\label{rq:twophi}
Consider two functions $\phi_1$ and $\phi_2$ which, up to an additive constant, coincide on some half-line: $\phi_1(s)=\phi_2(s)+\delta$ for all $s\geq r$ and some $\delta\in \re$. The optimal stopping problem (\ref{eq:optimal_stopping}) for $\phi_1$ has finite payoff if and only if the optimal stopping problem (\ref{eq:optimal_stopping}) for $\phi_2$ has finite payoff. Suppose this the case and denote respectively $\tau_1=\inf\{t\geq 0: X_t\leq g_*^1(S_t)\}$ and  $\tau_2=\inf\{t\geq 0: X_t\leq g_*^2(S_t)\}$ the stopping times, which yield these payoffs. Then $g_*^1(s)=g_*^2(s)$ for all $s\geq r$.
\end{remark}
We turn now to examine in more depth the relation between $g_*$ and  $\phi$. Consider an interval of constancy of $\phi$. Suppose that $\phi(\alpha)=\phi(\beta)$ for some $\alpha<\beta$ and $\phi$ increases on right of $\beta$ (or has a jump). Then starting the process at $x\in [\alpha,\beta)$ we face always the same problem: we get punished all the time, but get rewarded only if we attain the level $\beta$. Due to the Markov character of the diffusion, reasoning similarly to the derivation of Remark \ref{rq:formofthesolution}, we see that the optimal stopping time is then just an exit time of some interval $[\gamma,\beta]$ and $\gamma$ is then precisely the value of $g_*$ on $(\alpha,\beta)$.
\begin{remark}
\label{rq:constancy}
The functions $\phi$ and $g_*$, which induces the optimal stopping time, have the same intervals of constancy. In other words ${g_*'}_{|_{(\alpha,\beta)}}\equiv 0\Leftrightarrow\phi '_{|_{(\alpha,\beta)}}\equiv 0$.
\end{remark}
It is important to realize that if an interval of constancy of $\phi$ is long and the function doesn't grow fast enough afterwards, it may be optimal to stop immediately when coming to such a level of constancy. Thus the function $g_*$ may stay in part under and in part on the diagonal in $\re^2$.

Analysis of jumps of $\phi$ and $g_*$ is more involved. When $\phi$ is continuous, since the diffusion itself is continuous too, it's not hard to believe that $g_*$ should be continuous as well. It is maybe a little more tricky to convince yourself that if $\phi$ jumps, then $g_*$ has a jump too. We will try to argue this here, and moreover we will determine the height of the jump, following the reasoning of Meilijson developed in the proof of Corollary 1 in \cite{optimal_meilijson}. Suppose that $\phi$ has a jump in $s_0$: $\phi(s_0)-\phi(s_0-)=j>0$. For simplicity, suppose also that $\phi$ is constant on some interval before the jump: $\phi(s_0-\epsilon)=\phi(s_0-)$. Then, as noted above, starting the process at some $x\in [s_0-\epsilon,s_0)$ we wait till the first exit time of the interval $[g_*(x),s_0]$. If we exit at the bottom we stop, and if we exit at the top, we have a new, independent diffusion starting, for which the expected payoff is just $V(s_0,s_0)$.
For the process $X$ staring at $y\in [a,b]$, denote $\rho^y_{a,b}$ and $\rho^y_a$ respectively the first exit time of the interval $[a,b]$ and the first hitting time of the level $a$. We can then write the payoff $V(x,x)$ as
\begin{eqnarray}
  \label{eq:payoff_jump}
  V(x,x)&=&\e\Big[V(s_0,s_0)\mathbf{1}_{\{\rho^x_{g_*(x),s_0}=\rho^x_{s_0}\}}+\phi(x)\mathbf{1}_{\{\rho^x_{g_*(x),s_0}=\rho^x_{g_*(x)}\}}\nonumber\\&&-\int_0^{\rho^x_{g_*(x),s_0}}c(X_u)du\Big].
\end{eqnarray}
We know, by Remark \ref{rq:twophi}, that the value $g_*(s_0)$ is uniquely determined by the diffusion characteristics, the cost function $c$, and values of the function $\phi$ on $[s_0,\infty)$. Replace in the above display $g_*(x)$ by $a$ and denote this quantity $V^a(x,x)$. Then $V(x,x)=\sup_{a<s_0}V^a(x,x)$. This determines the value of $a_*^x=g_*(x)$ and therefore the height of the jump of $g_*$ at $s_0$. We will now show that the value of $a_*^x$ does not depend on $x$ (which is also clear by Remark \ref{rq:constancy}) and give an equation which describes it.

Actually the equation for $a^x_*$ is obvious: $a^x_*$ is the unique real, such that $\frac{\partial V^a(x,x)}{\partial a}_{|_{a=a^x_*}}=0$. This can be rewritten in the following equivalent manner, using the fact that $\p(\rho^x_{a,b}=\rho^x_b)=(L(x)-L(a))/(L(b)-L(a))$:
\begin{equation}
  \label{eq:diff_v_a}
L'(a_*^x)\times\frac{L(s_0)-L(x)}{(L(s_0)-L(a_*^x))^2}\Big(\phi(x)-V(s_0,s_0)\Big)=\frac{\partial}{\partial a}\e\Big[\int_0^{\rho^x_{a_*^x,s_0}}c(X_u)du\Big].
\end{equation}
We just have to show that actually $a_*^x$ does not depend on $x$. To this end let $a<y<x<s_0$ and write
\begin{eqnarray*}
  \e\Big[\int_0^{\rho^x_{a,s_0}}c(X_u)du\Big]&=&\e\Big[\int_0^{\rho^x_{s_0}}c(X_u)du\times \mathbf{1}_{\rho^x_{y,s_0}=\rho^x_{s_0}}\\
  &+& \int_0^{\rho^x_{y}}c(X_u)du\times \mathbf{1}_{\rho^x_{y,s_0}=\rho^x_{y}}\Big]\\
  &+& \frac{L(s)-L(x)}{L(s)-L(y)}\e\Big[\int_0^{\rho^y_{a,s_0}}c(X_u)du\Big],\nonumber
\end{eqnarray*}
where we used the strong Markov property at time $\rho^x_y$. The first two terms on the right-hand side do not depend on $a$, and differentiating the above equation with respect to $a$ yields an identity which proves that $a_*^x=a_*^y=:a_*$ for $x,y\in [s_0-\epsilon,s_0)$, since $\phi$ was supposed constant on this interval. It is also clear that $a_*$ depends only on the characteristics of $X$, the cost function $c$ and the payoff value $V(s_0,s_0)$. Since the diffusion is fixed, we write $a_*=a_*(c,V(s_0,s_0))$.

So far we only analyzed the dependence between $g_*$ and $\phi$, and it is the time to investigate the r\^ole of the cost function $c$. Suppose that the function $c$ disappears on some interval, $c_{|_{(\alpha,\beta)}}\equiv 0$. Thus wandering away from the maximum in the interval $(\alpha,\beta)$ doesn't cost us anything - there is no reason therefore to stop while in this interval. This implies that $(\alpha,\beta)$ is not in the range of $g_*$.
\begin{remark}
\label{rq:c}
  If $c_{|_{(\alpha,\beta)}}\equiv 0$ end $\tau_*$ yields the optimal payoff to the problem (\ref{eq:optimal_stopping}), then $(\alpha,\beta)\cap g_*(\re)=\emptyset$ and  $X_{\tau_*}\notin (\alpha,\beta)$ a.s.
\end{remark}
Thus, the impacts of the intervals of disappearance of $c$, and of the jumps of $\phi$, on the function $g_*$, are similar, only the former is much easier to describe. One could predict that the jumps of $c$ would provoke a similar behavior of $g_*$ as do the intervals of constancy of $\phi$. This however proves to be untrue. It stems from the fact that the values of $c$ are averaged through integrating them. Thus, jumps of $c$ will only produce discontinuities of $g_*'$. We will come back to this matter in Theorem \ref{thm:gen_max2}.

\section{Maximality principle revisited}
\label{sec:gen_max}

In this section we describe the main result of our paper, namely the complete solution to the optimal stopping problem (\ref{eq:optimal_stopping}). The solution is described in a sequence of three theorems with increasing generality of the form of the function $\phi$. In the first theorem we suppose that $\phi$ is of class $C^1$ and strictly increasing. Our theorem is basically a re-writing of the Theorem 3.1 found in Peskir \cite{MR2000c:60052}. Theorem \ref{thm:gen_max2} treats the case of a continuous function $\phi$, and is obtained form the previous one through an approximation procedure. Theorem \ref{thm:gen_max3} which deals with the general setup described in Section \ref{sec:notation}. It relays on Theorem \ref{thm:gen_max2} and on the work of Meilijson \cite{optimal_meilijson}, and is less explicit then Theorems \ref{thm:gen_max1} and \ref{thm:gen_max2}, as the treatment of jumps of $\phi$ is harder.
\begin{theorem}[Peskir]
\label{thm:gen_max1}
Let $\phi\in C^1$ with $\phi'>0$, and $c>0$ be continuous. The problem (\ref{eq:optimal_stopping}), has an optimal solution with finite payoff if and only if there exists a maximal solution $g_*$ of
\begin{equation}
  \label{eq:diff_g_general}
  g'(s)=\frac{\phi'(s)\sigma^2(g(s))L'(g(s))}{2c\big(g(s)\big)\big(L(s)-L(g(s))\big)},\quad a_X<s<b_X,
\end{equation}
  which stays strictly below the diagonal in $(a_X,b_X)^2$ (i.e. $g_*(s)<s$ for $a_X<s<b_X$). More precisely, let $s_*$ be such that $g_*(s_*)=a_X$. Then $g_*$ satisfies (\ref{eq:diff_g_general}) on $(s_*,b_X)$ and it is maximal such function (where the functions are compared on the interval, where both of them are superior to $a_X$), which stays below the diagonal.

In this case the payoff is given by
  \begin{equation}
    \label{eq:payoff}
    V_*(x,s)=\phi(s)+\int^{x}_{x\land g_*(s)}\big(L(x)-L(u)\big)c(u)m(du),
  \end{equation}
The stopping time $\tau_*=\inf\{t\geq 0: X_t\leq g_*(S_t)\}$ is then optimal whenever it satisfies (\ref{eq:optimal_stopping_cond}), otherwise it is ``approximately'' optimal\footnote{In the limit sense as in Peskir \cite{MR2000c:60052}.}.\\
Furthermore if there exists a solution $\rho$ of the optimal stopping problem (\ref{eq:optimal_stopping}) then $\p_{x,s}(\tau_*\leq \rho)=1$ and $\tau_*$ satisfies (\ref{eq:optimal_stopping_cond}).\\
If there is no maximal solution to (\ref{eq:diff_g_general}), which stays strictly below the diagonal in $\re^2$, then $V_*(x,s)=\infty$ for all $x\leq s$.\\
\end{theorem}
The last property for $\tau_*$ says it's pointwise the smallest solution, providing a uniqueness result. This remains true in more general setups and we will not repeat it below. This property follows from Peskir's \cite{MR2000c:60052} arguments but is also closely linked with properties of the Az\'ema-Yor stopping times, or any solutions to the Skorokhod embedding in general, we refer to Ob\l \'oj \cite{genealogia}, chapter 8, for details.
Note that for a general $\phi$, we defined in Section \ref{sec:notation} the function $\tilde{\phi}$, which coincides with $\phi$ on the interval $[r_\phi,\infty)$ and is of class $C^1$. In particular, we can apply the above theorem to the optimal stopping problem (\ref{eq:optimal_stopping}) with $\phi$ replaced by $\tilde{\phi}$. We denote $\tilde{g}_*$ the function which gives the optimal stopping time $\tilde{\tau}=\inf\{t\geq 0: X_t\leq \tilde{g}_*(S_t)\}$, which solves this problem.
\begin{theorem}
\label{thm:gen_max2}
  Let $\phi$ be as described in Section \ref{sec:notation}, but continuous, and $c>0$ be continuous. Then the optimal stopping problem (\ref{eq:optimal_stopping}) has a finite payoff if and only if the optimal stopping problem with $\phi$ replaced by $\tilde{\phi}$ has a finite payoff. In this case there exists a continuous function $g_*$, which satisfies the differential equation (\ref{eq:diff_g_general}) on the interior of the set $\re\setminus D_{\phi'}$, and coincides with $\tilde{g}_*$ on the interval $[r_\phi,\infty)$. The payoff is given by the formula (\ref{eq:payoff}) and is obtained for the stopping time $\tau=\inf\{t\geq 0: X_t\leq g_*(S_t)\}$, if it satisfies (\ref{eq:optimal_stopping_cond}). Otherwise this stopping time is approximately optimal.
\end{theorem}
It is important to note, that the function $g_*$ may not satisfy anymore $g_*(s)<s$, since if the constancy intervals of $\phi$ are too long it might be optimal to stop immediately (see also discussion around Remark \ref{rq:constancy}).

The formulation for the case of $\phi$ with discontinuities is somewhat more technical, as we have to define $g_*$ through an iteration procedure.
\begin{theorem}
\label{thm:gen_max3}
Let $\phi$ be as described in Section \ref{sec:notation}, and $c>0$ be continuous.
The optimal stopping problem (\ref{eq:optimal_stopping}) has finite payoff if and only if
the optimal stopping problem with $\phi$ replaced by $\tilde{\phi}$ has a finite payoff. In this case there exists a function $g_*$, continuous on the set $\re\setminus D_{\phi}$, which satisfies the following:
\begin{itemize}
\item $g_*(s)=\tilde{g}_*(s)$ for all $s\in [r_\phi,\infty)$;
\item $g_*$ is continuous on the interior of the set $\re\setminus D_\phi$ and differentiable on the interior of the set $\re\setminus D_{\phi'}$, where it satisfies the equation (\ref{eq:diff_g_general});
\item for all $s\in D_\phi$, $a_*=g_*(s-)$ satisfies the equation (\ref{eq:diff_v_a}) with $s_0=s$.
\end{itemize}
The payoff is described via (\ref{eq:payoff}).
\end{theorem}
It is important to note, that even though the optimal payoff $V$ appears in the equation (\ref{eq:diff_v_a}) the above construction of $g_*$ is feasible. Recall, that the jumps points were denoted $\ldots<d_1<d_0$. When determining the value of $g_*(d_i-)$ we already know the values of $g_*$ on $[d_i,\infty)$ and through the formula (\ref{eq:payoff}) also the payoff $V(x,s)$ for $d_i\leq x\leq s$.

\begin{proof} We now prove, in subsequent paragraphs, Theorems \ref{thm:gen_max1}, \ref{thm:gen_max2} and \ref{thm:gen_max3}.\\
\textit{Theorem \ref{thm:gen_max1}.}  \textit{Step 1.} We start with the case when the state space of the diffusion $X$ is the whole real line, $(a_X,b_X)=\re$, and the function $\phi$ is strictly increasing and of class $C^1$, and its image is the real line $\phi:\re\stackrel{on}{\to}\re$.
Let $Y_t=\phi(X_t)$. It is a diffusion with the scale function $L^Y(s)=L(\phi^{-1}(s))$, the drift coefficient $\delta^Y(y)=\phi'(\phi^{-1}(y))\delta(\phi^{-1}(y))+\frac{1}{2}\phi''(\phi^{-1}(y))\sigma^2(\phi^{-1}(y))$, the diffusion coefficient $\sigma^Y(y)=\phi'(\phi^{-1}(y))\sigma(\phi^{-1}(y))$ and the speed measure $m^Y(dy)=\frac{2dy}{L'(\phi^{-1}(y))\phi'(\phi^{-1}(y))\sigma^2(\phi^{-1}(y))}$. The state space of $Y$ is the whole real line.
We can rewrite the optimal stopping problem for $X$ in terms of $Y$ (where we add superscripts $X$ and $Y$ to denote the quantities corresponding to these two processes):
\begin{eqnarray}
  \label{eq:opt_x_opt_y}
  V^X_*(x,s)&=&\sup_{\tau}\e_{x,s}\Big(\phi(S^X_\tau)-\int_0^\tau c(X_t)dt\Big)\nonumber\\
  &=& \sup_{\tau}\e_{x,s}\Big(S^Y_\tau-\int_0^\tau c(\phi^{-1}(Y_t))dt\Big)\\
&=&V^Y_*(\phi(x),\phi(s))\textrm{ under }\tilde{\phi}=Id\textrm{ and }\tilde{c}=c\circ \phi^{-1}.\nonumber
\end{eqnarray}

We see therefore that the optimal solution for $X$ coincides with the optimal solution for $Y$ for a different problem. We can apply Theorem $3.1$ in Peskir \cite{MR2000c:60052} for this new problem for $Y$ to obtain that the optimal solution, if it exists, is given by
\begin{eqnarray}
  \label{eq:tau_x_tau_y}
  \tau_*&=&\inf\{t\geq 0: Y_t\leq g_*^Y(S^Y_t)\}\nonumber\\
  &=&\inf\Big\{t\geq 0: X_t\leq \phi^{-1}\Big(g_*^Y\big(\phi(S^X_t)\big)\Big)\Big\},
\end{eqnarray}
where $g^Y_*$ is the maximal solution which stays below the diagonal of the equation
\begin{equation}
\label{eq:diff_g_y}
  g_*^{Y'}(s)=\frac{\phi'\big(\phi^{-1}(g^Y_*(s))\big)^2\phi^{-1'}(g^Y_*(s))\sigma^2\big(\phi^{-1}(g_*^Y(s))\big)L'\big(\phi^{-1}(g^Y_*(s))\big)}{2c\big(\phi^{-1}(g^Y_*(s))\big)\Big(L(\phi^{-1}(s))-L\big(\phi^{-1}(g_*^Y(s))\big)\Big)}\;.
\end{equation}
We put $g_*^X(s)=\phi^{-1}\Big(g_*^Y\big(\phi(s)\big)\Big)$, so that the existence of $g^X_*$ is equivalent to the existence of $g^Y_*$.
We will now show that $g_*^X$ satisfies the equation (\ref{eq:diff_g_general}). We will use the following simple observations: $\phi^{-1'}(\phi(s))\phi'(s)=1$ and $\frac{1}{\phi^{-1'}(g^Y_*(\phi(s)))}=\phi'(g^X(s))$.
\begin{eqnarray}
  \label{eq:diff_g_x_diff_g_y}
  g_*^{X'}(s)&=&\phi^{-1'}\Big(g^Y_*(\phi(s))\Big)\times g^{Y'}_*(\phi(s))\times \phi'(s)\nonumber\\
  &=& \frac{\phi^{-1'}\Big(g^Y_*(\phi(s))\Big)^2\big(\phi'(g^X_*(s))\big)^2\sigma^2(g_*^X(s))L'(g^X_*(s))\phi'(s)}{2c\big(g_*^X(s)\big)\big(L(s)-L(g_*^X(s))\big)}\nonumber\\
  &=&\frac{\phi'(s)\sigma^2(g_*^X(s))L'(g^X_*(s))}{2c\big(g_*^X(s)\big)\big(L(s)-L(g_*^X(s))\big)}\;,\nonumber
\end{eqnarray}
in which we recognize (\ref{eq:diff_g_general}).
Note that $\big(g_*^X(s)<s$ for $s\in \re\big) \Leftrightarrow \big(g_*^Y(s)<s$ for $s\in \re$\big) as $\big(g_*^X(s)<s\big)\Leftrightarrow \big(g_*^Y(\phi(s))<\phi(s)\big)$ and $\phi$ is a one to one map from $\re$ to $\re$. This yields the description of $g^X_*$ as the maximal solution to (\ref{eq:diff_g_general}), which stays strictly below the diagonal in $\re^2$.

The expression for payoff follows also from Peskir \cite{MR2000c:60052} by a change of variables.

\noindent \textit{Step 2.} In this step we extend the previous results to the case, when the diffusion $X$ or the function $\phi$ are not unbounded. First, observe that the results of Peskir work just as well for a diffusion with an arbitrary state space $(a_X,b_X)$, $-\infty\leq a_X<b_X\leq \infty$. Indeed, if we take $\phi=Id$, the solution is given by the maximal solution to the equation (\ref{eq:diff_g_general}) for $s\in (a_X,b_X)$ and which stays strictly below the diagonal in $(a_X,b_X)^2$. Formally, we could repeat the proof in of Theorem 3.1 in \cite{MR2000c:60052}, in about the same manner that Peskir treats non-negative diffusions in Section 3.11 in the same article. Moreover, in this section, Peskir also discusses the boundary behavior, to which we will come below (notice that so far we considered unattainable boundaries, see Karlin and Taylor \cite[pp. 226--236]{MR82j:60003}).

We conclude therefore that the assumption on the state space of the diffusion, as well as on the range of $\phi$, were superficial. We need to keep $\phi'>0$, but we can have $\phi(\re)\subsetneq\re$. Note that $\phi(\re)$ is then an open interval $\phi(\re)=(a,b)$.
The equation (\ref{eq:diff_g_y}) is valid in $(a,b)$, but the original one (\ref{eq:diff_g_general}) is valid in $(a_X,b_X)$ as asserted, due to the fact that $g^X_*=\phi^{-1}\circ g^Y_*\circ \phi$.
The assumption, that $\phi$ is strictly increasing is needed for the new process to be a diffusion. We now turn to more general cases, which will be treated by approximation.
\smallskip

\noindent \textit{Theorem \ref{thm:gen_max2}.} We extend the previous results to functions $\phi$, which are just non-decreasing and continuous. This is done through an approximation procedure.\\
Let $\phi_n$ be an increasing sequence of strictly increasing $C^1$-functions, converging to $\phi$, such that $\phi_n'\to\phi'$, and $\phi_n|_{[r_\phi,\infty)}\equiv \phi|_{[r_\phi,\infty)}$. Note $V_n$ and $V$ the payoffs given in (\ref{eq:optimal_stopping}), under the condition (\ref{eq:optimal_stopping_cond}), for functions $\phi_n$ and $\phi$ respectively. Obviously $V_n\leq V_{n+1}\leq V$.  Furthermore, from Remark \ref{rq:twophi} it follows that the payoff $V$ is finite if and only if all $V_n$ are finite.
Suppose this is the case. We know, by Theorem \ref{thm:gen_max1} proved above, that  $V_n$ are given by a sequence of stopping times $\tau^n_*$, $\tau^n_*=\inf\{t\geq 0: X_t\leq g_*^n(S_t)\}$, where $g^n_*$ are defined as the maximal solutions of (\ref{eq:diff_g_general}), with $\phi$ replaced by $\phi_n$, which stay below the diagonal. It is easy to see that $\tau^n_*\leq \tau^{n+1}_*$ and therefore $g_*^n\geq g^{n+1}_*$. This is a direct consequence of the fact that in the problem (\ref{eq:optimal_stopping}) the reward $(\phi_n(S_\tau))$ increases and the punishment $(\int_0^\tau c(X_s)ds$) stays the same. It can also be seen from the differential equation (\ref{eq:diff_g_general}) itself.
We have therefore $\tau_*^n\nearrow \tau_*$ a.s., for some stopping time $\tau_*$, which is finite since the payoff associated with it is finite, $V<\infty$.
Moreover, $\tau_*=\inf\{t\geq 0:X_t\leq g_*(S_t)\}$, where $g^n_*\searrow g_*$. Note that the limit is finite and satisfies $g_*(s)<s$. Furthermore, we can pass to the limit in the differential equations describing $g^n_*$ to see that $g_*$ satisfies (\ref{eq:diff_g_general}).
Note that this agrees with Remark \ref{rq:constancy}.
Suppose that the stopping time $\tau_*$ satisfies (\ref{eq:optimal_stopping_cond}).
We then see in a straightforward manner, that $\tau_*$ solves the optimal stopping problem for $\phi$. It suffices to write:
\begin{eqnarray}
  \label{eq:limit_1}
  \e_{x,s}\Big(\phi(S_{\tau_*})-\int_0^{\tau_*} c(X_t)dt\Big)&\leq& V_*(x,s)\\&=&\sup_{\tau}\e_{x,s}\Big(\phi(S_\tau)-\int_0^\tau c(X_t)dt\Big)\nonumber \\
&\leq&\lim_{n\to \infty}\sup_{\tau}\e_{x,s}\Big(\phi_n(S_\tau)-\int_0^\tau c(X_t)dt\Big)\nonumber\\&=&\lim_{n\to\infty}\e_{x,s}\Big(\phi_n(S_{\tau_*^n})-\int_0^{\tau_*^n} c(X_t)dt\Big)\nonumber\\&=& \e_{x,s}\Big(\phi(S_{\tau_*})-\int_0^{\tau_*} c(X_t)dt\Big),\nonumber
\end{eqnarray}
where passing to the limit in both cases is justified by the monotone convergence of $\phi_n\nearrow \phi$ and $\tau^n_*\nearrow \tau_*$ a.s.. The expression for payoff given in (\ref{eq:payoff}) follows also upon taking the limit.

If the stopping time $\tau_*$ fails to satisfy (\ref{eq:optimal_stopping_cond}), we proceed as Peskir \cite[p. 1626]{MR2000c:60052}, to see that the expression (\ref{eq:payoff}) for payoff remains true (and we say $\tau^*$ is approximately optimal as it is in fact a limit of stopping times satisfying (\ref{eq:optimal_stopping_cond})).
\smallskip

\noindent\textit{Theorem \ref{thm:gen_max3}.} All that has to be verified is the expression for $g_*(d_i-)$, $i\geq 0$, where $D_\phi=\{d_0,d_1,\dots\}$.
Obviously, it is enough to discuss one jump, say $d_0$. The proposed value for $g(d_0-)$ is just the value we obtained at the end of Section \ref{sec:remarks}. We had however an assumption that $\phi$ was constant on some interval $[d_0-\epsilon,d_0)$ and we need to argue that it can be omitted. This is done by approximating $\phi$. Define $\phi_n$ through $\phi_n(s)=\phi(s)$ for all $s\in \re\setminus [d_0-\frac{1}{n},d_0)$ and $\phi_n$ constant and equal to $\phi(d_0-)$ on $[d_0-\frac{1}{n},d_0)$. The theorem then applies to functions $\phi_n$, as by Remark \ref{rq:twophi} the finiteness of the payoff for $\phi$, $\tilde{\phi}$ and $\phi_n$ are equivalent.
Passing to the limit, as in the proof of Theorem \ref{thm:gen_max2} above, yields the result. Seemingly, taking a sequence of continuous functions $\tilde{\phi}_n$ increasing to $\phi$ and repeating the reasoning in (\ref{eq:limit_1}) we see that the payoff equation (\ref{eq:payoff}) is still valid.
\end{proof}

Note that $g_*(b_X)=b_X$. It is evident when the state space is $(a_X,b_X]$, since we have to stop upon achieving the upper bound for the maximum $S$, as continuation will only decrease the payoff. Generally, if $g_*(b_X)=\beta<b_X$ (understood as limit if necessary), then as $S_t\to b_X$, the possible increase of payoff goes to zero but the cost till stopping doesn't.

The value of $g_*$ at $a_X$ depends on the boundary behavior of $X$ at $a_X$. Peskir  \cite{MR2000c:60052} dealt with this subject in detail. He considered the case $a_X=0$, as it was motivated by applications, but this has no impact on the result. We state briefly the results and refer to his work for more details. When $a_X$ is a natural or exit boundary point, then $g_*(a_X)=\lim_{s\downarrow a_X}g_*(s)=a_X$. In contrast, when $a_X$ is an entrance or regular instantaneously reflecting boundary then $g_*(a_X)=\lim_{s\downarrow a_X}g_*(s)<a_X$.

\section{General optimal Skorokhod embedding problem}
\label{sec:gen_skoro}

In this section we consider the so-called optimal Skorokhod embedding problem. The general idea is to design, for a given measure $\mu$, an optimal stopping problem with finite payoff, such that the stopping time which solves it embeds $\mu$, i.e.\ $X_\tau\sim \mu$.  We choose here to restrain ourselves to the classical setup of Brownian motion, $X=B$.
Study for general diffusions is equally possible, only the formulae would be more complicated and less-intuitive. In principle however it is not harder, since we know very well how to rephrase the classical Az\'ema-Yor embedding for the setup of any real diffusion (Az\'ema and Yor \cite{MR82c:60073a}, Cox and Hobson \cite{cox_hobson}, see Ob\l \'oj \cite{genealogia} for a complete description of the subject).

This problem was first introduced in the context of classification of contingent claims, where the measure $\mu$ had an interpretation in terms of risk associated with an option (see Peskir \cite{MR1982093}).
Peskir in his article on the optimal Skorokhod embedding \cite{MR2001e:60083} gave an explicit solution, but he assumed for simplicity that the measure $\mu$ had positive density on $\re$. One could think that this will generalize easily to arbitrary measures. Unfortunately this proves not to be true. In the setup when $\phi(x)=x$, it is not possible to solve the optimal embedding problem for a measure, which has an atom in the interior of its support. To see this, note that the presence of an atom in the interior of the support implies that the barycentre function $\Psi_\mu$ has a jump, which in turn means that the function $g_*=\Psi_\mu^{-1}$ has to be constant on some interval. However, the derivative of $g$ is given by $g'(s)=[2c(g(s))(s-g(s))]^{-1}$ which is strictly positive.
In this sense, it is the presence of the function $\phi$ which proves fundamental to solve the problem. Meilijson \cite{optimal_meilijson}, who had $c(x)=c$ but arbitrary $\phi$, solved this problem for any centered measure $\mu$ with finite variance, but his solution is not really explicit.
We will see that actually to obtain a general explicit solution, both function $c$ and $\phi$ are needed, as the former allows to treat the regular (absolutely continuous) part of a measure, and the latter serves to obtain atoms.

We set ourselves two goals. Firstly, for a given measure $\mu$, we want to identify all the pairs $(\phi,c)$ such that the optimal stopping time $\tau_*$ which solves (\ref{eq:optimal_stopping}) embeds $\mu$ in $B$, $B_{\tau_*}\sim\mu$. We are interested in stopping times $\tau_*$ such that $(B_{t\land \tau_*}:t\ge 0)$ is a uniformly integrable martingale. This limits the class of measures which are admissible to $\mu$ such that $\int_{\re}|x|d\mu(x)<\infty$ and $\int_{\re}xd\mu(x)=0$.
Our second goal is to give a particular explicit description of a particular pair $(\phi,c)$.

Let us start with the case of a measure $\mu$ with a positive density $f>0$ on $\re$. Then, we know that for $\tau_*=\inf\{t\geq 0: B_t\leq g_*(S_t)\}$ we have $B_{\tau_*}\sim \mu$ if and only if $g_*(s)=\Psi_\mu^{-1}(s)$, where $\Psi_\mu^{-1}$ denotes the inverse of $\Psi_\mu$. This is a direct consequence of the Az\'ema-Yor embedding \cite{MR82c:60073a} and the fact that $g$ is increasing with $g(s)<s$.
Let us investigate the differential equation satisfied by $\Psi_\mu^{-1}$.
Using the explicit formula (\ref{eq:ay_def}), we obtain
\begin{equation}
  \label{eq:diff_psi_inv}
  \Big(\Psi_\mu^{-1}(s)\Big)'=\frac{\overline{\mu}(\Psi_\mu^{-1}(s))}{(s-\Psi_\mu^{-1}(s))f\big(\Psi_\mu^{-1}(s)\big)}.
\end{equation}
Comparing this with the differential equation for $g_*$ given by (\ref{eq:diff_g_general}), we see that we need to have
\begin{equation}
  \label{eq:cond_pair}
  \frac{\phi'(s)}{2c(\Psi_\mu^{-1}(s))}=\frac{\overline{\mu}(\Psi_\mu^{-1}(s))}{f(\Psi_\mu^{-1}(s))}.
\end{equation}
Equivalently, we have
\begin{equation}
  \label{eq:cond_pair}
  \frac{\phi'(\Psi_\mu(u))}{2c(u)}=\frac{\overline{\mu}(u)}{f(u)},\quad \textrm{for }u\in \re.
\end{equation}
Finally, we have to ensure that the payoff is finite. This is easy since we want $\e \phi(S_{\tau_*})$ to be finite and we know that $S_{\tau_*}=\Psi_\mu(B_{\tau_*})$, which implies that we need to have
\begin{equation}
  \label{eq:cond_pair_phi}
\int \phi(\Psi_\mu(x))d\mu(x)<\infty.
\end{equation}
We have thus obtained the following proposition:
\begin{prop}
  \label{prop:opt_skoro}
Let $\mu$ be a probability measure on $\re$ with a strictly positive density $f$.
Then the optimal stopping problem (\ref{eq:optimal_stopping}) has finite payoff and the stopping time $\tau_*$, which yields it embeds $\mu$ in $B$, i.e. $B_{\tau_*}\sim \mu$ if and only if $(\phi,c)$ satisfy (\ref{eq:cond_pair}) and (\ref{eq:cond_pair_phi}) for all $s\geq 0$. The stopping time $\tau_*$ is then just the Az\'ema-Yor stopping time given by (\ref{eq:ay_st_def}).
\end{prop}
Note that, if we take $\phi(s)=s$ so that $\phi'(s)=1$ we obtain a half of the hazard function for $c$, which is the result of Peskir, stated in Theorem \ref{thm:peskir2}. The integrability condition (\ref{eq:cond_pair_phi}) then
reads $\int \Psi_\mu(x)d\mu(x)<\infty$, which is known (see Az\'ema and Yor \cite{MR82c:60073b}) to be equivalent to the $L\log L$ integrability condition on $\mu$: $\int_1^\infty x\log x d\mu(x)<\infty$.
\begin{proof}
  In light of the reasoning the led to (\ref{eq:cond_pair}) all we have to comment on is the continuity of the function $c$. Theorems \ref{thm:gen_max1}-\ref{thm:gen_max3} were formulated for a continuous function $c>0$. Here, keeping the condition $c>0$, we drop the assumption on continuity. Still, this is not a problem. It suffices to take a sequence of functions $c_n\searrow c$ a.e.\ and proceed as in the proof of Theorem \ref{thm:gen_max2} around (\ref{eq:limit_1}), to see that the function $g_*$, which yields the solution $\tau_*=\inf\{t\geq 0: B_t\leq g_*(S_t)\}$ to the optimal stopping problem (\ref{eq:optimal_stopping}), satisfies locally the equation (\ref{eq:diff_g_general}) and is continuous. This in turn implies that it is indeed the inverse of the barycentre function associated with $\mu$.
\end{proof}

Identifying all pairs $(\phi,c)$, which solve the optimal Skorokhod embedding problem for an arbitrary measure is harder. More precisely it's just not explicit any more. The reason is hidden in Theorem \ref{thm:gen_max3} - the description of jumps of $g_*$, which correspond to intervals not charged by the target measure, is done through an iteration procedure and is not explicit. This is exactly the reason why Meilijson was only able to prove the existence of the solution to the optimal Skorokhod embedding without giving explicit formulae. In our approach we will use the duality between $\phi$ and $c$ to encode explicitly both the jumps and the regular part of the target measure. Still, we will not be able to cover all probability measures.

Let $\mu$ be a centered probability measure $\int_{\re} |x|d\mu(x)<\infty$, $\int_{\re}x d\mu(x)=0$ and note $-\infty\leq a_\mu<b_\mu\leq +\infty$ respectively the lower, and the upper, bound of its support. Suppose $\mu$ is a sum of its regular and atomic parts: $\mu=\mu_r+\mu_a$, $d\mu_r(x)=f(x)dx$ and $\mu_a=\sum_{i\in \zr}p_i \delta_{j_i}$. In other words $f$ is the density of the absolutely continuous (with respect to the Lebesgue measure) part and
$\dots<j_{-1}<j_0<j_1<\ldots$ are the jump points of $\mu$, which are also the jump points of the barycentre function $\Psi_\mu$, so that $\Psi_\mu(\re)=\re_+\setminus \Big(\bigcup_{i\in \zr}\big(\Psi_\mu(j_i),\Psi_\mu(j_i+)\big]\Big)$. We may note that $\Psi\mu(j_i+)=\Psi_\mu(j_i)+\frac{p_i(\Psi_\mu(j_i)-j_i)}{\mu((j_i,+\infty))}$.
\begin{theorem}
\label{thm:opt_skoro}
  In the above setup, in the case when $f>0$ or $f\geq 0$ but $\mu_a=0$, the optimal stopping problem (\ref{eq:optimal_stopping}) with
$$c(x)=\left \{
    \begin{array}{ll}
      \frac{f(x)}{2\overline{\mu}(x)}&\textrm{, for }x\in[a_\mu,b_\mu]\\
      +\infty&\textrm{, for }x\in \re\setminus [a_\mu,b_\mu]
    \end{array}\right.\quad\textrm{and}\quad
\left \{
    \begin{array}{l}
\phi'(x)=\mathbf{1}_{\Psi_\mu(\re)}(x),\\
\phi\textrm{ continuous, }\phi(0)=0
    \end{array}\right.
$$
has finite payoff under (\ref{eq:cond_pair_phi}). The payoff is realized by the Az\'ema-Yor stopping time $\tau_*=\inf\{t\geq 0: S_t\geq \Psi_\mu(B_t)\}$ and $B_{\tau_*}\sim\mu$.
\end{theorem}
Note that the condition $\phi(0)=0$ is just a convention, as adding a constant to $\phi$ doesn't affect the solution $\tau_*$ of (\ref{eq:optimal_stopping}) (cf.\ Remark \ref{rq:twophi}). The restriction (\ref{eq:cond_pair_phi}), thanks to the definition of $\phi$, is satisfied in particular when $\int \Psi_\mu(x)d\mu(x)<\infty$, that is if $\mu$ satisfies the $L\log L$ integrability condition.
\begin{proof}
  The problem arising form discontinuity of $c$ is treated exactly as in the proof of Proposition \ref{prop:opt_skoro} above. In the case of $f>0$ on $\re$ the above theorem follows immediately from Theorem \ref{thm:gen_max2}. If $\mu$ has no atoms, then $\Psi_\mu(\re)=\re_+$ and so $\phi=Id$. From the differential equation (\ref{eq:diff_g_general}) we can derive a differential equation for the inverse of $g_*$, from which it is clear that we can treat intervals where $c$ is zero by passing to the limit.
Thus, our theorem is also valid in the case of absolutely continuous measure but with a density which can be zero at some intervals.
Finally, the case of atoms with $\limsup_{n\to \infty} j_n =+\infty$ is also handled upon approximating and taking the limit. Indeed this does not bother us here as we do not need anymore to describe the limiting solution through a differential equation, as was the case for Theorems \ref{thm:gen_max1} - \ref{thm:gen_max3}.
\end{proof}
Unfortunately, our method does not work for an arbitrary measure $\mu$. If one tries to apply it for a purely atomic measure, he would end up with $\phi\equiv c\equiv 0$. In fact, when there is no absolutely continuous part we need to take $\phi$ discontinuous, as Meilijson does. We are not able to solve the problem explicitly then. We cannot treat measures with singular non-atomic component.

\section{Important inequalities}
\label{sec:ieq}

In the above section we saw how the solution to the optimal stopping problem yields a solution to the so-called optimal Skorokhod embedding problem. This actually motivated our research but is by no means a canonical application of the maximality principle. Probably the main and the most important applications are found among stopping inequalities. We will try to present some of them here. A specialist in the optimal stopping theory will find nothing new in this section and can probably skip it, yet we think it is useful to put it in this note, as it completes our study and allows us to give some references for further research.
The main interest of the method presented here, is that the constants obtained are always optimal. Some of the inequalities below, as (\ref{eq:ineq4}), are easy to obtain with 'some' constant. It was however the question of the optimal constant which stimulated researchers for a certain time.

Consider a continuous local martingale $X=(X_t:t\geq 0)$ and a stopping time $T$ in the natural filtration of $X$, $T<\infty$ a.s. We have then
\begin{equation} \label{eq:ineq1} \e S^X_T\leq \sqrt{\e X_T^2},\quad\textrm{where }S_t^X=\sup_{s\leq t}X_s,
\end{equation} and this is optimal. This simple inequality was first observed by Dubins and Schwarz \cite{MR89m:60101}. It can be easily seen in the following manner (localizing and passing to the limit if necessary):  $\e S_T=\e [S_T-X_T]\leq \sqrt{\e[S_T-X_T]^2}=\sqrt{\e[X_T^2]}$, where we used the fact that $(S_t-X_t)^2-X_t^2$ is a local martingale (see Ob\l\'oj and Yor \cite{ja_yor_maxmart} for various applications of these martingales). This inequality is optimal as the equality is attained for $T=\inf\{t\geq 0: S^X_t-X_t=a\}$ for any $a>0$ ($X_T$ has then shifted exponential distribution with parameter $\frac{1}{a}$).
If we want to establish an analogous inequality for $|X|$, that is for a submartingale, some more care is needed. We propose to follow Peskir \cite{MR99i:60157, MR2000c:60052} in order to obtain more general inequalities.

Consider the optimal stopping problem (\ref{eq:optimal_stopping}), with $X=|B|$ the absolute value of a Brownian motion, $\phi(s)=s^p$ and $c(x)=cx^{p-2}$, where $p>1$. We solve the problem applying Theorem \ref{thm:gen_max1}. The differential equation (\ref{eq:diff_g_general}) takes the form
\begin{equation}
\label{eq:ineq1.5}
g'(s)=\frac{ps^{p-1}}{2c g(s)^{p-2}(s-g(s))},
\end{equation}
which is solved by $g(s)=\alpha s$, where $\alpha$ is the larger root of the equation $\alpha^{p-1}-\alpha^p=\frac{p}{2c}$, which is seen to have solutions for $c\geq p^{p+1}/2(p-1)^{(p-1)}$. Taking $c\searrow p^{p+1}/2(p-1)^{(p-1)}$ yields
\begin{equation}
\label{eq:ineq2}
\e\Big(\max_{0\leq t\leq \tau} |B_\tau|^p\Big)\leq \Big(\frac{p}{p-1}\Big)^{p}\e|B_\tau|^p,
\end{equation}
where $\tau$ is any stopping time in the natural filtration of $|B|$ such that $\e \tau^{p/2}<\infty$. Note that to obtain the right-hand side of the above inequality we used the fact that $\e\big(\int_0^\tau |B_t|^{p-2}\big)=\frac{2}{p(p-1)}\e |B_\tau|^{p-1}$. More generally we could consider the optimal stopping problem (\ref{eq:optimal_stopping}) with $\phi(s)=s^p$ and $c(x)=cx^{q-1}$ for $0<p<q+1$, $q>0$. Optimizing upon $c$ Peskir \cite{MR99i:60157} obtains
\begin{equation}
\label{eq:ineq3}
\e\Big(\max_{0\leq t\leq \tau} |B_\tau|^p\Big)\leq \gamma^*_{p,q}\e\Big(\int_0^\tau |B_t|^{q-1}dt\Big)^{p/(q-1)},
\end{equation}
for all stopping times $\tau$ of $|B|$. The optimal constant $\gamma^*_{p,q}$ is given implicitly as a solution to an equation. In the case $p=1$ it can be written explicitly, $\gamma^*_{1,q}=\big(q(1+q)/2\big)^{1/(1+q)}\big(\Gamma(2+1/q)\big)^{q/(1+q)}$. This was obtained independently by Jacka \cite{MR89j:60054} and Gilat \cite{MR89m:60102}.
In particular we find $\gamma^*_{1,1}=\sqrt{2}$ which is exactly the value found by Dubins and Schwarz \cite{MR89m:60101}. If we consider the well-known bounds
\begin{equation}
\label{eq:ineq4}
c_p \e\Big(\int_0^\tau |B_t|^{p-2}dt\Big)\leq \e\Big(\max_{0\leq t\leq \tau} |B_\tau|^p\Big)\leq C_p \e\Big(\int_0^\tau |B_t|^{p-2}dt\Big),
\end{equation}
where $p>1$ and $c_p$, $C_p$ are some universal constants, we see that the inequality (\ref{eq:ineq3}) complements them for $0<p\leq 1$, where (\ref{eq:ineq4}) wouldn't have much sense. Also, as pointed out above, the method presented here allows to recover the optimal constants.

The inequalities given above are just a sample of applications of the method presented in this note. Peskir \cite{MR2000c:60052}, for example, develops also a $L \log L$-type inequalities and some inequalities for Geometric Brownian motion. Although the method presented is very general it has of course its limits. For example it seems one could not recover optimal constants\footnote{These can be recovered using different methodology involving determinist time changing (see Pedersen and Peskir \cite{MR2001i:60076}).} $a_p$ and $A_p$, obtained by Davis \cite{MR54:6260}, such that for any stopping time $T$
\begin{eqnarray}
  \label{eq:ineq5}
&&\e|B_T|^p\leq A_p\e T^{p/2},\textrm{ for } 0<p<\infty,\textrm{ and}\\
 && a_p\e T^{p/2}\leq \e |B_t|^p,\textrm{ for }  1<p<\infty,\; \e T^{p/2}<\infty .
\end{eqnarray}
The maximality principle deals with the maximum process and is an essential tool in its study. However the passage from the maximum to the terminal value is quite complicated and in general we should not hope that optimal inequalities for the terminal value could be obtained from the ones for maximum.
\smallskip

{\small\noindent \textbf{Acknowledgment.} Author wants to express his gratitude towards professors Goran Peskir and Marc Yor for their help and support.
\def\cprime{$'$} \def\polhk#1{\setbox0=\hbox{#1}{\ooalign{\hidewidth
  \lower1.5ex\hbox{`}\hidewidth\crcr\unhbox0}}}

\end{document}